\documentclass{amsproc}%
\usepackage{graphicx}
\usepackage{amscd}
\usepackage{amsmath}
\usepackage{amsfonts}
\usepackage{amssymb}%
\theoremstyle{plain}

\numberwithin{equation}{section}

\begin{document}
\title[Linear algebra in lattices ]{Linear algebra in lattices, the Fitting lemma}
\author{Jen\H{o} Szigeti}
\address{Institute of Mathematics, University of Miskolc, Miskolc, Hungary 3515}
\email{jeno.szigeti@uni-miskolc.hu}
\thanks{\noindent This work was supported by OTKA of Hungary No. T043034 and K61007}

\begin{abstract}
Lattice theoretical generalizations of some classical linear algebra results
are formulated. A vector space is replaced by its subspace lattice and a
linear map is replaced by the induced lattice map. This map is a complete join
homomorphism and satisfies some additional natural conditions. The object of
our study is a pair consisting of a complete lattice and one of its complete
join endomorphisms. We prove the Fitting lemma for such a pair.

\end{abstract}
\maketitle

\noindent1. PRELIMINARIES

\bigskip

\noindent The purpose of this paper is to develop linear algebra for a pair
$(L,\lambda)$ satisfying some of the so called JNB conditions, here
$(L,\vee,\wedge,0,1)$ is a complete lattice and $\lambda:L\longrightarrow L$
is a complete $\vee$-homomorphism with $\lambda(0)=0$. A lattice theoretical
Jordan normal base theorem for a nilpotent $\lambda$ was proved in [3]. In the
present work we continue the investigations started in [3], to provide a self
contained treatment, we repeat the proof of a lemma which we need from [3].

\noindent Let $R$ be a ring and $\varphi:M\longrightarrow M$ an $R$%
-endomorphism of the unitary left $R$-module $_{R}M$ (a vector space together
with a linear map also can be considered). For a submodule $N\leq M$ let
$\varphi(N)=\{\varphi(x)\mid x\in N\}$ denote the $\varphi$-image of $N$.
Clearly, $\varphi(N)\leq M$ is an $R$-submodule and for any family $N_{\gamma
}$, $\gamma\in\Gamma$ of $R$-submodules in $M$ we have%
\[
\varphi\left(  \underset{\gamma\in\Gamma}{\sum}N_{\gamma}\right)
=\underset{\gamma\in\Gamma}{\sum}\varphi(N_{\gamma}).
\]
If $X,Y\leq M$ are $R$-submodules with $X\subseteq Y$ and $\varphi
(X)=\varphi(Y)$, then%
\[
X+(Y\cap\ker(\varphi))=Y
\]
and $\varphi(Y\cap\ker(\varphi))=\{0\}$ holds for the $R$-submodule $Y\cap
\ker(\varphi)\leq M$.

\noindent If $X\leq M$ is an $R$-submodule with $X\subseteq\varphi(N)$, then
$\varphi^{-1}(X)=\{u\in M\mid\varphi(u)\in X\}$ and we have $\varphi(Z)=X$ for
the $R$-submodule $Z=N\cap\varphi^{-1}(X)\subseteq N$ of $M$.

\noindent The above observations were used in [3] to formulate the so called
JNB2 and JNB3 properties of the map $\lambda$\ (or of the pair $(L,\lambda)$),
where
\[
L=\text{Sub}(_{R}M)
\]
is the complete (and modular) lattice of the $R$-submodules of $M$ with
respect to the containment relation $\subseteq$\ and%
\[
\lambda(N)=\varphi(N)
\]
is the induced $L\longrightarrow L$ complete $\vee$-homomorphism (we note that
$\underset{\gamma\in\Gamma}{\vee}N_{\gamma}=\underset{\gamma\in\Gamma}{\sum
}N_{\gamma}$).

\bigskip

\noindent2. FINITE\ ATOMIC BASE IN A LATTICE

\bigskip

\noindent Let $(L,\vee,\wedge,0,1)$\ be a lattice with $\leq_{L}$ ($<_{L}$)
being the induced (strict) partial order on $L$. A finite set $\{a_{i}%
\mid1\leq i\leq n\}$ of atoms in $L$ is called \textit{independent}, if%
\[
a_{\pi(1)}<_{L}a_{\pi(1)}\vee a_{\pi(2)}<_{L}...<_{L}a_{\pi(1)}\vee a_{\pi
(2)}\vee...\vee a_{\pi(n)}%
\]
is a strictly ascending chain for some permutation $\pi$ of $\{1,2,...,n\}$.
An \textit{atomic base} of $L$ is an independent set of atoms $\{a_{i}%
\mid1\leq i\leq n\}$ with $a_{1}\vee a_{2}\vee...\vee a_{n}=1$.

\noindent If $1\in L$ is a join of atoms and $1=b_{1}\vee b_{2}\vee...\vee
b_{m}$, ($b_{1},b_{2},...,b_{m}$ are atoms in $L$) with $m$ being the smallest
possible, then $b_{1}\vee b_{2}\vee...\vee b_{m}$ is an irredundant join and%
\[
b_{\sigma(1)}<_{L}b_{\sigma(1)}\vee b_{\sigma(2)}<_{L}...<_{L}b_{\sigma
(1)}\vee b_{\sigma(2)}\vee...\vee b_{\sigma(m)}%
\]
is strictly ascending for any choice of the permutation $\sigma$ of
$\{1,2,...,m\}$. Thus $\{b_{i}\mid1\leq i\leq m\}$ is an atomic base of $L$.

\noindent The following trivial statement about the join of atoms is already
present in [3].

\bigskip

\noindent\textbf{2.1.Proposition.} \textit{Let }$(L,\vee,\wedge,0)$%
\textit{\ be a lattice with the following properties}.

\noindent(i)\textit{ The cover relation }$x\vartriangleleft x\vee a$\textit{
holds for all }$x\in L$\textit{ and for all atoms }$a\in L$\textit{ with
}$x\neq x\vee a$\textit{ (atomic cover property, it is satisfied in upper
semimodular lattices).}

\noindent(ii)\textit{ Any two composition series of the form}%
\[
0=x_{0}\vartriangleleft x_{1}\vartriangleleft...\vartriangleleft
x_{n}\text{\textit{ and }}0=y_{0}\vartriangleleft y_{1}\vartriangleleft
...\vartriangleleft y_{m}%
\]
\textit{with }$x_{n}=y_{m}$\textit{ are of the same length: }$n=m$\textit{.}

\noindent\textit{If }$\{a_{i}\mid1\leq i\leq n\}$\textit{ is an independent
set of atoms in }$L$\textit{, then }$a_{1}\vee a_{2}\vee...\vee a_{n}$\textit{
is an irredundant join (or equivalently}%
\[
a_{1}\vee a_{2}\vee...\vee a_{n}=a_{1}\oplus a_{2}\oplus...\oplus a_{n}%
\]
\textit{is a direct sum).}

\bigskip

\noindent Let $(L,\vee,\wedge,0,1)$\ be a complete lattice and $\lambda
:L\longrightarrow L$ a complete $\vee$-homomorphism with $\lambda(0)=0$. The
\textit{image} and the \textit{kernel} of $\lambda$ can be defined as follows:%
\[
w=\text{im}(\lambda)=\lambda(1)\text{ and }z=\ker(\lambda)=\vee\{x\in
L\mid\lambda(x)=0\}.
\]
Clearly, the $\vee$-property of $\lambda$ ensures that $\lambda(x)\leq
_{L}\lambda(1)=w$ for all $x\in L$ and that $\lambda(x)=0\Longleftrightarrow
x\leq_{L}z$.

\noindent Now we repeat the formulation of the properties (mentioned in the
introduction) of the induced map $N\longmapsto\varphi(N)$ in the language of
lattices. We keep the notations JNB2 and JNB3 of [3].

\bigskip

\noindent JNB2: For any choice of the elements $x\leq_{L}y$ with
$\lambda(x)=\lambda(y)$ we can find an element $u\in L$ such that $y=x\vee u$
and $\lambda(u)=0$.

\bigskip

\noindent JNB3: For any $x\in L$ the map $\lambda:[0,x]\longrightarrow
\lbrack0,\lambda(x)]$ is surjective.

\bigskip

\noindent\textbf{2.2.Lemma.}\textit{ If every element of }$L$\textit{ is a
join of atoms and }$\lambda$\textit{\ satisfies JNB3, then any atom }%
$a\in\lbrack0,\lambda(1)]$\textit{ can be written as }$a=\lambda(a^{\ast}%
)$\textit{, where }$a^{\ast}\in L$\textit{ is an atom.}

\bigskip

\noindent\textbf{Proof.} Since $\lambda:[0,1]\longrightarrow\lbrack
0,\lambda(1)]$ is surjective, we can find an element $x\in L$ with
$a=\lambda(x)$. Now we have $x=\vee\{b_{i}\mid i\in I\}$, where each $b_{i}\in
L$ is an atom. Clearly, $b_{i}\leq_{L}x$ implies that $0\leq_{L}\lambda
(b_{i})\leq_{L}\lambda(x)=a$. For an index $i\in I$, we have either
$\lambda(b_{i})=0$ or $\lambda(b_{i})=a$. If $\lambda(b_{i})=0$ holds for all
$i\in I$, then we obtain%
\[
a=\lambda(x)=\lambda(\underset{i\in I}{\vee}b_{i})=\underset{i\in I}{\vee
}\lambda(b_{i})=0,
\]
a contradiction. Thus $\lambda(b_{j})=a$ holds for some $j\in I$.$\square$

\bigskip

\noindent\textbf{2.3.Proposition.} \textit{Let }$(L,\vee,\wedge,0,1)$%
\ \textit{be a complete lattice and }$\lambda:L\longrightarrow L$\textit{ a
complete }$\vee$\textit{-homomorphism such that }$\lambda(0)=0$.\textit{
Consider the interval sublattices }$[0,w]$\textit{ and }$[0,z]$\textit{\ of
}$L$\textit{ with }$w=$im$(\lambda)$\textit{ and }$z=\ker(\lambda)$\textit{.}

\noindent\textit{\noindent If }$\{b_{j}\mid1\leq j\leq q\}$\textit{ is
independent set of atoms in }$[0,z]$\textit{ and }$\{a_{i}\mid1\leq i\leq
p\}$\textit{ is a set of atoms in }$L$\textit{ such that }$\{\lambda
(a_{i})\mid1\leq i\leq p\}$\textit{ is an independent set of atoms (in
}$[0,w]$\textit{), then}%
\[
\{b_{j}\mid1\leq j\leq q\}\cup\{a_{i}\mid1\leq i\leq p\}
\]
\textit{is independent in }$L$\textit{.}

\noindent\textit{\noindent If }$\lambda$\textit{\ satisfies JNB2, }$b_{1}\vee
b_{2}\vee...\vee b_{q}=z$\textit{ and }$\lambda(a_{1})\vee\lambda(a_{2}%
)\vee...\vee\lambda(a_{p})=\lambda(1)$\textit{, then}%
\[
b_{1}\vee b_{2}\vee...\vee b_{q}\vee a_{1}\vee a_{2}\vee...\vee a_{p}=1,
\]
\textit{i.e. }$\{b_{j}\mid1\leq j\leq q\}\cup\{a_{i}\mid1\leq i\leq
p\}$\textit{ is an atomic base of }$L$\textit{.}

\bigskip

\noindent\textbf{Proof.} Now we have strictly ascending chains%
\[
b_{\alpha(1)}<_{L}b_{\alpha(1)}\vee b_{\alpha(2)}<_{L}...<_{L}b_{\alpha
(1)}\vee b_{\alpha(2)}\vee...\vee b_{\alpha(q)}=b
\]
in $[0,z]$ and%
\[
\lambda(a_{\beta(1)})<_{L}\lambda(a_{\beta(1)})\vee\lambda(a_{\beta(2)}%
)<_{L}...<_{L}\lambda(a_{\beta(1)})\vee\lambda(a_{\beta(2)})\vee...\vee
\lambda(a_{\beta(p)})
\]
in $[0,w]$. We claim that%
\[
b<_{L}b\vee a_{\beta(1)}<_{L}b\vee a_{\beta(1)}\vee a_{\beta(2)}<_{L}%
...<_{L}b\vee a_{\beta(1)}\vee a_{\beta(2)}\vee...\vee a_{\beta(p)}%
\]
is also strictly ascending. Indeed $1\leq i\leq p$ and%
\[
b\vee a_{\beta(1)}\vee a_{\beta(2)}\vee...\vee a_{\beta(i-1)}=b\vee
a_{\beta(1)}\vee a_{\beta(2)}\vee...\vee a_{\beta(i)}%
\]
would imply that%
\[
a_{\beta(i)}\leq_{L}b\vee a_{\beta(1)}\vee a_{\beta(2)}\vee...\vee
a_{\beta(i-1)}.
\]
Since $b\in\lbrack0,z]$ ensures that $\lambda(b)=0$, the $\vee$-property of
$\lambda$ gives that%
\[
\lambda(a_{\beta(i)})\leq_{L}\lambda(b)\vee\lambda(a_{\beta(1)})\vee
\lambda(a_{\beta(2)})\vee...\vee\lambda(a_{\beta(i-1)})=
\]%
\[
=\lambda(a_{\beta(1)})\vee\lambda(a_{\beta(2)})\vee...\vee\lambda
(a_{\beta(i-1)}),
\]
whence%
\[
\lambda(a_{\beta(1)})\vee...\vee\lambda(a_{\beta(i-1)})=\lambda(a_{\beta
(1)})\vee...\vee\lambda(a_{\beta(i-1)})\vee\lambda(a_{\beta(i)})
\]
follows, a contradiction.

\noindent If $\lambda(a_{1})\vee\lambda(a_{2})\vee...\vee\lambda
(a_{p})=\lambda(1)$, then we have%
\[
\lambda(a_{1}\vee a_{2}\vee...\vee a_{p})=\lambda(1)
\]
and JNB2 implies that%
\[
(a_{1}\vee a_{2}\vee...\vee a_{p})\vee u=1
\]
holds for some $u\in L$ with $\lambda(u)=0$. Since $u\leq_{L}z$, we have%
\[
b_{1}\vee b_{2}\vee...\vee b_{q}\vee a_{1}\vee a_{2}\vee...\vee a_{p}%
=z\vee(a_{1}\vee a_{2}\vee...\vee a_{p})=1
\]
proving that $\{b_{j}\mid1\leq j\leq q\}\cup\{a_{i}\mid1\leq i\leq p\}$ is an
atomic base of $L$.$\square$

\bigskip

\noindent\textbf{2.4.Theorem.} \textit{Let }$(L,\vee,\wedge,0,1)$\textit{\ be
a complete lattice and }$\lambda:L\longrightarrow L$\textit{ a complete }%
$\vee$\textit{-homomorphism with }$\lambda(0)=0$\textit{ such that }$\lambda
$\textit{\ satisfies JNB2 and JNB3. Assume that every element of }%
$L$\textit{\ is a join of atoms and consider the interval sublattices }%
$[0,w]$\textit{ and }$[0,z]$\textit{\ of }$L$\textit{ with }$w=$im$(\lambda
)$\textit{ and }$z=\ker(\lambda)$\textit{. If }$[0,w]$\textit{ has a }%
$p$\textit{ element atomic base and }$[0,z]$\textit{ has a }$q$\textit{
element atomic base, then we can find a }$p+q$\textit{ element atomic base in
}$L$\textit{.}

\bigskip

\noindent\textbf{Proof.} If $\{c_{i}\mid1\leq i\leq p\}$ is a $p$ element
atomic base of $[0,w]$, then Lemma 2.2 ensures the existence of atoms
$a_{i}\in L$, $1\leq i\leq p$ such that $\lambda(a_{i})=c_{i}$ for all $1\leq
i\leq p$. Take a $q$ element atomic base $\{b_{j}\mid1\leq j\leq q\}$\ of
$[0,z]$, then $\{b_{j}\mid1\leq j\leq q\}\cup\{a_{i}\mid1\leq i\leq p\}$ is a
$p+q$ element atomic base of $L$ by Proposition 2.3.$\square$

\bigskip

\noindent\textbf{Remark.} Theorem 2.4 is the generalization of the well known%
\[
\dim(\text{im}\varphi)+\dim(\ker\varphi)=\dim V
\]
property of the linear map $\varphi:V\longrightarrow V$.

\bigskip

\noindent3. THE\ FITTING\ LEMMA

\bigskip

\noindent\textbf{3.1.Proposition.} \textit{Let }$(L,\vee,\wedge,0,1)$%
\textit{\ be a complete lattice and }$\lambda_{1},\lambda_{2}:L\longrightarrow
L$\textit{ complete }$\vee$\textit{-homomorphisms with }$\lambda
_{1}(0)=\lambda_{2}(0)=0$\textit{. Then }$\lambda_{1}\circ\lambda
_{2}:L\longrightarrow L$\textit{ is a complete }$\vee$\textit{-homomorphisms
with }$(\lambda_{1}\circ\lambda_{2})(0)=0$\textit{ and we have the following
implications.}

\begin{enumerate}
\item \textit{If }$\lambda_{1}$\textit{ and }$\lambda_{2}$\textit{ satisfy
JNB3, then }$\lambda_{1}\circ\lambda_{2}$\textit{\ satisfies JNB3.}

\item \textit{If }$\lambda_{1}$\textit{ satisfies JNB2, }$\lambda_{2}$\textit{
satisfies JNB2 and JNB3, then }$\lambda_{1}\circ\lambda_{2}$\textit{ satisfies
JNB2.}

\item \textit{If }$\lambda_{1}$\textit{ and }$\lambda_{2}$\textit{ both
satisfy JNB2 and JNB3, then }$\lambda_{1}\circ\lambda_{2}$\textit{\ satisfies
JNB2 and JNB3.}
\end{enumerate}

\bigskip

\noindent\textbf{Proof.}

\noindent1. If $x\in L$, then $\lambda_{2}:[0,x]\longrightarrow\lbrack
0,\lambda_{2}(x)]$ and $\lambda_{1}:[0,\lambda_{2}(x)]\longrightarrow
\lbrack0,\lambda_{1}(\lambda_{2}(x))]$ are surjective maps, thus their
composition $\lambda_{1}\circ\lambda_{2}:[0,x]\longrightarrow\lbrack
0,\lambda_{1}(\lambda_{2}(x))]$ is also surjective.

\noindent2. If $x\leq_{L}y$ and $\lambda_{1}(\lambda_{2}(x))=\lambda
_{1}(\lambda_{2}(y))$ holds for $x,y\in L$, then $\lambda_{2}(x)\leq
_{L}\lambda_{2}(y)$ and JNB2 of $\lambda_{1}$ gives an element $u\in L$ such
that $\lambda_{2}(y)=\lambda_{2}(x)\vee u$ and $\lambda_{1}(u)=0$. Now
$u\in\lbrack0,\lambda_{2}(y)]$ and JNB3 of $\lambda_{2}$\ gives an element
$v\in\lbrack0,y]$ such that $\lambda_{2}(v)=u$. In view of $x\vee v\leq_{L}y$
and%
\[
\lambda_{2}(x\vee v)=\lambda_{2}(x)\vee\lambda_{2}(v)=\lambda_{2}(x)\vee
u=\lambda_{2}(y)
\]
JNB2 of $\lambda_{2}$\ gives an element $w\in L$ such that $(x\vee v)\vee w=y$
and $\lambda_{2}(w)=0$. Thus $x\vee(v\vee w)=y$ and%
\[
\lambda_{1}(\lambda_{2}(v\vee w))=\lambda_{1}(\lambda_{2}(v)\vee\lambda
_{2}(w))=\lambda_{1}(u\vee0)=\lambda_{1}(u)=0
\]
prove the JNB2 property of $\lambda_{1}\circ\lambda_{2}$.$\square$

\bigskip

\noindent\textbf{3.2.Corollary.} \textit{Let }$(L,\vee,\wedge,0,1)$%
\textit{\ be a complete lattice and }$\lambda:L\longrightarrow L$\textit{ a
complete }$\vee$\textit{-homomorphism with }$\lambda(0)=0$\textit{. If
}$\lambda$\textit{ satisfies JNB2 and JNB3, then any power }$\lambda^{k}%
$\textit{ satisfies JNB2 and JNB3.}

\bigskip

\noindent\textbf{3.3.Proposition.} \textit{Let }$(L,\vee,\wedge,0,1)$%
\textit{\ be a complete lattice and }$\lambda:L\longrightarrow L$\textit{ a
complete }$\vee$\textit{-homomorphism with }$\lambda(0)=0$\textit{. We have
the following implications.}

\begin{enumerate}
\item im$(\lambda^{k})=$im$(\lambda^{k+1})$\textit{ implies that }%
im$(\lambda^{i})=$im$(\lambda^{k})$\textit{ for all integers }$i\geq
k$\textit{.}

\item $\ker(\lambda^{k})=\ker(\lambda^{k+1})$\textit{ implies that }%
$\ker(\lambda^{i})=\ker(\lambda^{k})$\textit{ for all integers }$i\geq
k$\textit{.}

\item \textit{If }$\lambda$\textit{ satisfies JNB2, then }im$(\lambda
)=$im$(\lambda^{2})$\textit{ implies that }im$(\lambda)\vee\ker(\lambda
)=1$\textit{.}

\item \textit{If }$\lambda$\textit{ satisfies JNB3, then }$\ker(\lambda
)=\ker(\lambda^{2})$\textit{ implies that }im$(\lambda)\wedge\ker(\lambda
)=0$\textit{.}
\end{enumerate}

\bigskip

\noindent\textbf{Proof.}

\noindent1. Clearly, $\lambda^{k}(1)=\lambda^{k+1}(1)$ implies that
$\lambda^{k+1}(1)=\lambda^{k+2}(1)$.

\noindent2. Take $x=\ker(\lambda^{k+2})$, then $\lambda^{k+1}(\lambda
(x))=\lambda^{k+2}(x)=0$ implies $\lambda(x)\leq_{L}\ker(\lambda^{k+1})$,
whence $\lambda(x)\leq_{L}\ker(\lambda^{k})$ can be obtained. Thus $x\leq
_{L}\ker(\lambda^{k+1})$ is a consequence of $\lambda^{k+1}(x)=\lambda
^{k}(\lambda(x))=0$. Since $\ker(\lambda^{k+1})\leq_{L}\ker(\lambda^{k+2})$ is
obvious, we get $\ker(\lambda^{k+1})=\ker(\lambda^{k+2})$.

\noindent3. Using $w=$im$(\lambda)=\lambda(1)$ and $z=\ker(\lambda)$, we have%
\[
\lambda(w\vee z)=\lambda(w)\vee\lambda(z)=\lambda(\lambda(1))\vee0=\lambda
^{2}(1)=\lambda(1).
\]
The JNB2 property of $\lambda$ gives that $(w\vee z)\vee u=1$ and
$\lambda(u)=0$ for some $u\in L$. Thus $u\leq_{L}z$ and $w\vee z=(w\vee z)\vee
u=1$.

\noindent4. Take $y=$im$(\lambda)\wedge\ker(\lambda)$, then $y\leq_{L}%
\lambda(1)$ and thus $y\in\lbrack0,\lambda(1)]$. The JNB3 property of
$\lambda$ ensures the existence of an element $v\in L$ such that
$\lambda(v)=y$. Now $y\leq_{L}\ker(\lambda)$ implies that $\lambda
^{2}(v)=\lambda(y)=0$, whence $v\leq_{L}\ker(\lambda^{2})$ can be derived.
Since $\ker(\lambda^{2})=\ker(\lambda)$, we have $y=\lambda(v)=0$.$\square$

\bigskip

\noindent\textbf{3.4.Theorem.} \textit{Let }$(L,\vee,\wedge,0,1)$\textit{\ be
a complete lattice and }$\lambda:L\longrightarrow L$\textit{ a complete }%
$\vee$\textit{-homomorphism with }$\lambda(0)=0$\textit{ such that }$\lambda
$\textit{\ satisfies JNB2 and JNB3. If }$L$\textit{\ is both Artinian and
Noetherian (or }$L$\textit{\ is of finite height), then there exists an
integer }$r\geq1$\textit{ such that }im$(\lambda^{r})\oplus\ker(\lambda
^{r})=1$\textit{ (i.e. }im$(\lambda^{r})\vee\ker(\lambda^{r})=1$\textit{ and
}im$(\lambda^{r})\wedge\ker(\lambda^{r})=0$\textit{).}

\bigskip

\noindent\textbf{Proof.} The Artinian condition ensures that%
\[
\text{im}(\lambda)\geq_{L}\text{im}(\lambda^{2})\geq_{L}...\geq_{L}%
\text{im}(\lambda^{i})\geq_{L}...
\]
is not a strictly descending chain, thus im$(\lambda^{k})=$im$(\lambda^{k+1})$
holds for some integer $k\geq1$.

\noindent The Noetherian condition ensures that%
\[
\ker(\lambda)\leq_{L}\ker(\lambda^{2})\leq_{L}...\leq_{L}\ker(\lambda^{i}%
)\leq_{L}...
\]
is not a strictly ascending chain, thus $\ker(\lambda^{l})=\ker(\lambda
^{l+1})$ holds for some integer $l\geq1$.

\noindent Take $r=\max\{k,l\}$, then%
\[
\text{im}(\lambda^{r})=\text{im}(\lambda^{2r})\text{\ and }\ker(\lambda
^{r})=\ker(\lambda^{2r})
\]
by part 1 and part 2 of Proposition 3.3, respectively. In view of Corollary
3.2, the power $\lambda^{r}$ satisfies JNB2 and JNB3. Since $\lambda
^{2r}=(\lambda^{r})^{2}$, the application of part 3 and part 4 of Proposition
3.3 give that im$(\lambda^{r})\vee\ker(\lambda^{r})=1$ and im$(\lambda
^{r})\wedge\ker(\lambda^{r})=0$.$\square$

\bigskip

\noindent\textbf{Remark.} Theorem 3.4 is the generalization of the well known
Fitting lemma. We note that a completely different lattice theoretical Fitting
lemma can be found in [2], it can be considered as an analogue and not a
generalization of the classical result.

\bigskip

\noindent REFERENCES

\bigskip

\begin{enumerate}
\item Anderson, F.W. and Fuller, K.R.: \textit{Rings and Categories of
Modules}, GTM Vol. 13, Springer Verlag, New York, 1974.

\item K\"{o}rtesi, P. and Szigeti, J.:\textit{ A general approach to the
Fitting lemma}, Mathematika (London) Vol. 52 (2005), 155-160.

\item Szigeti, J.: \textit{The Jordan normal base in lattices and nilpotent
endomorphisms of finitely generated semisimple modules,} submitted
\end{enumerate}

\end{document}